\input amstex
\loadeufm

\documentstyle{amsppt}

\magnification=\magstep1

\baselineskip=20pt
\parskip=5.5pt
\hsize=6.5truein
\vsize=9truein
\NoBlackBoxes

\define\br{{\Bbb R}}

\define\e{{\varepsilon}}
\define\OO{{\Omega}}
\define\CL{{\Cal{L}}}
\define\bu{{\bold{u}}}
\define\bbf{{\bold{f}}}
\define\bo{{\bold{0}}}
\define\bv{{\bold{v}}}
\define\bg{{\bold{g}}}
\define\bw{{\bold{w}}}
\topmatter
\title
The $L^p$ Dirichlet Problem for Elliptic Systems
 on Lipschitz Domains
\endtitle

\author Zhongwei Shen
\endauthor

\leftheadtext{Zhongwei Shen}
\rightheadtext{Dirichlet Problem}
\address Department of Mathematics, University of Kentucky,
Lexington, KY 40506.
\endaddress

\email shenz\@ms.uky.edu
\endemail

\abstract
We develop a new approach to the $L^p$ Dirichlet problem
via $L^2$ estimates and reverse H\"older inequalities.
We apply this approach to second order elliptic systems
and the polyharmonic equation on 
 a bounded Lipschitz domain $\OO$ in $\br^n$.
For $n\ge 4$ and
$2-\e< p<\frac{2(n-1)}{n-3} +\e$,
we establish the solvability of the Dirichlet problem
 with boundary data in
$L^p(\partial\OO)$.
In the case of the polyharmonic equation $\Delta^\ell u=0$
with $\ell\ge 2$,
the range of $p$ is sharp if $4\le n\le 2\ell +1$.
\endabstract



\endtopmatter

\document

\centerline{\bf 1. Introduction}

Let $\OO$ be a bounded Lipschitz domain in $\br^n$. Consider
the system of second order elliptic operators
$(\CL(\bu))^r =-a_{ij}^{rs} D_iD_j u^s$,
where $D_i=\partial/\partial x_i$ and
$r,s=1,\dots,m$.
We assume that the coefficients $a_{ij}^{rs}$,
$1\le i,j\le n$, $1\le r,s \le m$ are real constants
satisfying the symmetry condition $a_{ij}^{rs}
=a_{ji}^{sr}$ and the Legendre-Hadamard ellipticity condition:
$$
\mu|\xi|^2|\eta|^2
\le a_{ij}^{rs} \xi_i\xi_j \eta^r\eta^s
\le
\frac{1}{\mu} |\xi|^2|\eta|^2,
\tag 1.1
$$
for some $\mu>0$ and any $\xi\in \br^n$, $\eta\in\br^m$.
The primary purpose of this paper is to study the Dirichlet problem
$$
\left\{
\aligned
& \CL(\bu)=\bo\  \text{ in }\OO,\\
& \bu=\bbf \in L^p(\partial\OO)
 \text{ on } \partial\OO\ \text{ and }
\   (\bu)^*\in L^p(\partial\OO),
\endaligned
\right.
\tag 1.2
$$
where $(\bu)^*$ denotes the nontangential maximal function of
$\bu$, and the boundary value of $\bu$ is taken in the sense of
nontangential convergence.
The following is one of main results of the paper.

\proclaim{Theorem 1.3}
Let $\OO$ be a bounded Lipschitz domain in $\br^n$, $n\ge 4$ with
connected boundary. Then there exists $\e>0$ depending only on 
$n$, $m$, $\mu$ and $\OO$ such that, given any
$\bbf \in L^p(\partial\OO)$ with $2-\e
< p< \frac{2(n-1)}{n-3} +\e$,
 the Dirichlet
problem (1.2) has a unique solution. Moreover, the solution
$\bu$ satisfies the estimate $\| (\bu)^*\|_p\le C\, \| \bbf\|_p$.
\endproclaim

With Theorem 1.3 at our disposal, using area integral estimates for
elliptic systems \cite{DKPV} and a duality argument found in \cite{V1},
we also establish the solvability of the regularity problem
with boundary data in  $L^p_1(\partial\OO)$ for $p$ in the dual range.
Here $L^p_1(\partial\OO)$ is the space of functions in $L^p(\partial\OO)$
whose first order (tangential) derivatives are also in $L^p$.

\proclaim{\bf Theorem 1.4}
Let $\OO$ be a bounded Lipschitz domain in $\br^n$, $n\ge 4$ with
connected boundary. Then there exists $\e_1>0$ depending
only on $n$, $m$, $\mu$ and $\OO$ such that, given
any $\bbf\in L^p_1(\partial\OO)$ with
$\frac{2(n-1)}{n+1}-\e_1 <p< 2+\e_1$, there exists a unique $\bu$
satisfying $\CL(\bu)=\bo$ in $\OO$, $\bu=\bbf$ on $
\partial\OO$, and
$(\nabla \bu)^*\in L^p(\partial\OO)$.
Moreover, the solution $\bu$ may
be represented by a single layer potential, and
we have $\| (\nabla \bu)^*\|_p\le C\, \| \nabla_t \bbf\|_p$,
where $\nabla_t \bbf$ denotes the tangential derivatives
of $\bbf$ on $\partial\OO$.
\endproclaim

We remark that for Laplace's equation in Lipschitz domains, 
the Dirichlet and Neumann problems
with boundary data in $L^p$ are well understood.
Indeed, it has been known since the early 1980's that the 
Dirichlet problem with optimal estimate $ (u)^*\in L^p$
is uniquely solvable for $2-\e<p\le \infty$, and the Neumann problem
as well as the regularity problem
with optimal estimate $(\nabla u)^*\in L^p$ is uniquely solvable
 for $1<p<2+\e$
(see \cite{D1, D2, JK, V1, DK1, K1}).
For elliptic systems as well as higher order 
elliptic equations, the solvability of the $L^p$ boundary value
problems were established for $n\ge 3$ and $2-\e <p< 2+\e$
(see \cite{DKV1, DKV2, FKV, F, V3, G, PV3}). 
This was achieved by the method of layer potentials.
The main tool was certain Rellich-Payne-Weinberg-Necas 
identities. 
We mention that for Laplace's equation, the $L^p$ estimate 
($2<p<\infty$) for
the Dirichlet problem follows directly from the $L^2$ estimate and the
well known maximal principle, by interpolation.
The $L^p$ estimate ($1<p<2$) for the Neumann problem relies
on the classical H\"older estimates for solutions of second order
elliptic equations of divergence form with bounded measurable
coefficients.
Similar $L^\infty$ and H\"older estimates, however,
 are not readily available
for elliptic systems or higher-order elliptic equations
on nonsmooth domains.

Nevertheless, in the case $n=3$, the $L^p$ boundary value problems
for the optimal ranges of $p$ were solved for elliptic systems
\cite{DK2, S1, S2} and higher order elliptic equations \cite{PV1, PV2, PV4} 
(also see \cite{MM} for systems on manifold).
In particular, it was proved that the Dirichlet problem (1.2)
is uniquely solvable for $2-\e< p\le \infty$.
The basic idea in \cite{PV2} is to estimate the decay rate of the
Green's function on $\OO$, using Rellich 
identities, Caccioppoli inequalities as well as the fact that the 
$L^p$ Dirichlet
problem is solvable for some $p<2$.
If $n=3$ (or $2$), the decay rate obtained is fast enough to yield the
Miranda-Agmon maximum principle $\|u\|_{L^\infty (\OO)}
\le C\, \| u\|_{L^\infty(\partial\OO)}$, from which the
$L^p$ estimates follows.
Recently
in \cite {S3, S4},  we  extended this approach to higher dimensions. 
However, instead of $L^p$ estimates, we were only able to 
establish dilation-invariant estimates in terms of the nontangential
maximal functions in the Morrey spaces
and their pre-duals, and certain weighted $L^2$ estimates with power
weights.

In this paper we develop a new approach to the $L^p$
Dirichlet problem. The basic idea is as
follows. First we prove a reverse H\"older inequality
on $\partial\OO$:
$$
\left(\frac{1}{|\Delta_r|}\int_{\Delta_r}
|(\bu)^*|^p d\sigma\right)^{1/p}
\le C\, 
\left(\frac{1}{|\Delta_{5r}|}
\int_{\Delta_{5r}} |(\bu)^*|^2 d\sigma\right)^{1/2},
\tag 1.5
$$
 for $
2<p<\frac{2(n-1)}{n-3} +\e,
$
where $\Delta_r$ is a surface ``cube'' on $\partial\OO$, and
$\bu$ satisfies $\CL \bu =\bo$ in $\OO$,
$(\bu)^*\in L^2(\partial\OO)$ and $\bu =\bo$
on $\Delta_{6r}$.
To do this, we use the inequality $(\bu)^*\le \sup_K |\bu|
+C\, I_1((\nabla \bu)^*)$, where $K$ is a compact set in $\OO$ and
$I_1$ denotes the fractional integral of order one on $\partial\OO$.
The restriction of $p$ comes from the fact that $(\nabla \bu)^*$
is in $L^q$ only for some $q>2$. Inequality (1.5)
should be considered as a localization estimate
in the sense that if the boundary values of two solutions
$\bu$, $\bold{v}$ agree on $\Delta_{6r}$, then
$(\bu-\bold{v})^*$ behaves far better than either $(\bu)^*$ 
or $(\bold{v})^*$
on $\Delta_r$.

Let $M_{\Delta_{2r}}$ denote the Hardy-Littlewood maximal function
on $\partial\OO$, localized to the surface cube
$\Delta_{2r}$. The second step in the proof of Theorem 1.3,  which is
motivated
by a paper of Caffarelli and Peral \cite{CP},  is to
establish a good-$\lambda$ type inequality by
using a real variable argument and estimate (1.5).
Indeed, 
for $\lambda\ge \lambda_0$ and the unique 
 solution of the $L^2$ Dirichlet
problem $\CL \bu =\bo$ in $\OO$ and $\bu=\bbf$ on $\partial\OO$,
we will show that
$$
|E(A\lambda)|\le \delta |E(\lambda)|
+|\left\{ P\in \Delta_r: \ 
M_{\Delta_{2r}}(|\bbf|^2)(P)>\gamma \lambda\right\}|,
\tag 1.6
$$
where $E(\lambda)=\left\{ P\in \Delta_r : \ M_{\Delta_{2r}}(|(\bu)^*|^2)(P)
>\lambda\right\}$, and $A$, $\delta$ and $\gamma$ are positive constants
with the property $\delta A^{p/2}<1$. 
The desired estimate $\|(\bu)^*\|_p\le C\, \| \bbf\|_p$
follows readily from
(1.6) by an integration in $\lambda$.

The proof of estimate (1.6) relies on a Calder\'on-Zygmund
decomposition. The key observation is that
if 
$$
\frac{1}{|\Delta_{12r}|}
\int_{\Delta_{12r}}
|\bu|^2 \, d\sigma\le \gamma \lambda
\ \text{ and }\
\frac{1}{|\Delta_{12r}|}
\int_{\Delta_{12r}}
|(\bu)^*|^2\, d\sigma\le \lambda,
\tag 1.7
$$
then the localization estimate (1.5) implies that
$$
|\left\{ P\in \Delta_r:\ 
M_{\Delta_{2r}}\big(|(\bu)^*|^2\big)(P)> A\lambda\right\}|
\le
C\, \left(\frac{1}{A^{p/2}}+\frac{\gamma}{A}\right) \, |\Delta_r|.
\tag 1.8
$$
At the end of section 2, we use these techniques to formulate
a theorem of independent interest
 on the $L^p$ boundedness of operators on $\br^n$. 
It essentially states that if a sublinear operator $T$ is bounded
on $L^2(\br^n)$, and satisfies a $L^p$ localization 
property in the spirit of (1.5) for some $p>2$, then it is bounded on
$L^q(\br^n)$ for all $2<q<p$
(see Theorem 2.32).
This theorem may be considered as a refined version of the
well-known Calder\'on-Zygmund Lemma
on weak type $(1,1)$ bounds of  
$L^2$ bounded operators.

We point out that the main ingredients of this new approach
are, (1) interior estimates for solutions of $\CL \bu =0$,
(2) the solvability of the $L^2$ Dirichlet problem,
(3) the solvability of the regularity problem
with data in $L^q_1(\partial\OO)$ for some $q>2$.
Thus our method also applies to the higher order elliptic equations
studied in \cite{DKV1, V3, PV3}, where 
 the Dirichlet and regularity problems in $L^p$  are solved
for $2-\e<p<2+\e$.
As an example, we consider the Dirichlet problem
for the polyharmonic equation:
$$
\left\{
\aligned
&\Delta^\ell u =0 \ \ \text{ in } \OO,\\
&D^\alpha u =f_\alpha \ \text{for } |\alpha|\le \ell-2,\ 
\frac{\partial^{\ell-1} u}{\partial N^{\ell-1}} 
=g\  \text{ on }\partial\OO,\  
 (\nabla^{\ell-1} u)^*\in L^p(\partial\OO),
\endaligned
\right.
\tag 1.9
$$
where $\ell\ge 2$. In (1.9), $\alpha=(\alpha_1,
\dots,\alpha_n)$ is a multi-index, $D^\alpha
=D_1^{\alpha_1}D_2^{\alpha_2}\cdots D_n^{\alpha_n}$,
and $|\alpha|=\alpha_1+\cdots +\alpha_n$.
Also $N$ denotes the outward unit normal to
$\OO$, and $\frac{\partial^{\ell-1}u}{\partial N^{\ell-1}}
=\sum_{|\alpha|=\ell-1}
\frac{(\ell-1)!}{\alpha!} N^\alpha D^\alpha u$.

\proclaim{Theorem 1.10}
Let $\OO$ be a bounded Lipschitz  domain in $\br^n$, $n\ge 4$
with connected boundary.
Then there exists $\e>0$ depending on $n$, $\ell$ and $\OO$
such that if $f=\{ f_\alpha:\, 0\le |\alpha|\le \ell-2\}
\in WA^p_{\ell-1}(\partial\OO)$ and
$g\in L^p(\partial\OO)$ with $2-\e < p<\frac{2(n-1)}{n-3}+\e$, the
Dirichlet problem (1.9) has a unique solution.
Moreover, the solution $u$ satisfies
$$
\| (\nabla^{\ell-1} u)^*\|_p\le C\, \left\{
\| g\|_p +\sum_{|\alpha|= \ell -2}
\|\nabla_t f_\alpha\|_p \right\}.
\tag 1.11
$$
\endproclaim

We remark that $WA^p_{\ell-1}(\partial\OO)$ is a Sobolev
space of the so-called Whitney arrays
 on $\partial\OO$. We refer the reader to \cite{V3} or \cite{PV3}
for its definition.
As we mentioned earlier, in the case $n=3$,
the $L^p$ Dirichlet problem
(1.9) was solved in \cite{PV1, PV2, PV4} for the optimal
range $2-\e< p\le \infty$.
It was also pointed out in \cite{PV3} that 
estimate (1.11) fails 
 in general for $p>\frac{2(n-1)}{n-3}$ and $4\le n\le 2\ell+1$, 
or $p>\frac{2\ell}{\ell-1}$ and $n\ge 2\ell+2$.
Thus the range for $p$ in Theorem 1.10
is sharp in the case $4\le n\le 2\ell+1$.

The paper is organized as follows. Theorem 1.3 is proved in section 
2. 
The proof of Theorem 1.10, which is very similar to
that of Theorem 1.3, is sketched in section 3.
Finally we give the proof of Theorem 1.4 in section 4.

\medskip

\centerline{\bf 2. The Dirichlet Problem}

Let $\OO$ be a bounded Lipschitz domain in $\br^n$.
Assume $0\in \partial\OO$ and 
$$
\OO\cap B(0,r_0)
=\left\{ (x^\prime, x_n)\in \br^n:\ \
x_n>\psi(x^\prime)\right\} \cap B(0,r_0),
\tag 2.1
$$
where $B(0,r_0)$ denotes the ball centered at $0$ with
radius $r_0$, $\psi:\br^{n-1}\to\br$ is a Lipschitz function.
For $r>0$,  we let
$$
\aligned
&\Delta_r =\big\{ (x^\prime, \psi(x^\prime))
\in\br^n:\ \ 
|x_1|<r,\dots, |x_{n-1}|<r\big\},\\
&D_r =\big\{ (x^\prime,x_n)\in \br^n:\ \ 
|x_1|<r,\dots, |x_{n-1}|<r, \ \psi(x^\prime)<x_n< \psi(x^\prime) +r\big\}.
\endaligned
\tag 2.2
$$
Note that if $0<r<c\, r_0$, $\Delta_r\subset\partial\OO$ and
$D_r\subset \OO$.

\proclaim{\bf Lemma 2.3} Let $0<r<c\, r_0$.
There exists $\e>0$ depending only on $n$, $m$, $\mu$ and $\OO$,
such that if 
 $\CL \bu =\bo$ in $\OO$, $(\bu)^*\in L^2(\partial\OO)$ and
$\bu =\bo$ on $\Delta_{8r}$, then $(\bu)^*\in L^p(\Delta_r)$ and 
$$
\left(\frac{1}{|\Delta_r|}
\int_{\Delta_r}
|(\bu)^*|^p d\sigma\right)^{1/p}
\le C\, 
\left(\frac{1}{|\Delta_{5r}|}
\int_{\Delta_{5r}}
|(\bu)^*|^2 d\sigma\right)^{1/2},
\tag 2.4
$$
where $2<p<2(n-1)/(n-3)+\e$ and $|\Delta|$ denotes the surface 
measure of $\Delta$.
\endproclaim

\demo{Proof}
Recall that the nontangential maximal function of $\bu$ is defined by
$$
(\bu)^* (P)
=\sup \big\{ |\bu(x)|:\ \ x\in \OO \text { and }\ x\in \gamma(P)\big\}\ \ 
\text{ for } P\in \partial\OO,
\tag 2.5
$$
where $\gamma (P)=\big\{ x\in \OO: \ \ |x-P|<2\, \text{dist}(x,\partial\OO)
\big\}$. Let
$$
\aligned
&\Cal{M}_1 (\bu)(P)
=\sup\big\{ |\bu(x)|:\ \ x\in \gamma(P) \ \ \text{ and }
\ \ |x-P|\le c\, r\big\},\\
&
\Cal{M}_2 (\bu)(P)
=\sup\big\{ |\bu(x)|:\ \ x\in \gamma(P) \ \ \text{ and }
\ \ |x-P|> c\, r\big\}.
\endaligned
\tag 2.6
$$
Then $(\bu)^*=\max \big\{ \Cal{M}_1 (\bu), \Cal{M}_2(\bu) \big\}$.
We first estimate $\Cal{M}_2(\bu)$. To do this,
 we use interior estimates
to obtain
$$
|\bu(x)|\le \frac{C}{r^n}
\int_{B(x,c r)} |\bu (y)| dy
\le \frac{C}{r^{n-1}}
\int_{\Delta_{5r}} |(\bu)^*| d\sigma,
\tag 2.7
$$
where $x\in \gamma (P)$, $P\in \Delta_r$ and $|x-P|\ge cr$.
It follows that for any $p>2$,
$$
\aligned
&\left(\frac{1}{|\Delta_r|}
\int_{\Delta_r} |\Cal{M}_2(\bu)|^p d\sigma\right)^{1/p}
\le \sup_{\Delta_r}
\Cal{M}_2(\bu)\\
&\le \frac{C}{|\Delta_{5r}|}
\int_{\Delta_{5r}}
|(\bu)^*| d\sigma
\le 
C\, \left(\frac{1}{|\Delta_{5r}|}
\int_{\Delta_{5r}}
|(\bu)^*|^2 d\sigma\right)^{1/2}.
\endaligned
\tag 2.8
$$

Next, to estimate $\Cal{M}_1(\bu)$ on $\Delta_r$, we write
$$
\bu(x^\prime, x_n)-\bu(x^\prime, \widetilde{{x}}_n)
=-\int_{x_n}^{\widetilde{x}_n} \frac{\partial \bu}{\partial s}
(x^\prime, s) \, ds.
\tag 2.9
$$
Let $K= \big\{ (x^\prime, x_n)\in \OO:\ \ 
|x^\prime|<c\, r,\ \psi(x^\prime) +c\, r<x_n
<\widetilde{c}\, r \big\}$.
Using (2.9) and 
interior estimates for $\nabla \bu$, it is not hard to show that
for $P\in \Delta_r$,
$$
\Cal{M}_1(\bu)(P)
\le \sup_K |\bu|
+C\, \int_{\Delta_{2r}}
\frac{ (\nabla \bu)^*_{D_{2r}}(Q)}{|P-Q|^{n-2}}
\, d\sigma(Q),
\tag 2.10
$$
where $(\nabla \bu)^*_{D_{2r}}$ denotes the nontangential maximal function
of $\nabla \bu$ with respect to the Lipschitz sub-domain $D_{2r}$.
By the fractional integral estimate, this implies that
$$
\left(\frac{1}{|\Delta_r|}
\int_{\Delta_r}
|\Cal{M}_1(\bu)|^p d\sigma\right)^{1/p}
\le \sup_K |\bu|
+C\, r
\left(\frac{1}{|\Delta_{2r}|}
\int_{\Delta_{2r}}|(\nabla \bu)^*_{D_{2r}}|^q d\sigma\right)^{1/q},
\tag 2.11
$$
where $\frac{1}{p}=\frac{1}{q}-\frac{1}{n-1}$ and $1<q<n-1$.
We now choose $q>2$, depending only on $n$, $m$, $\mu$ and $\OO$, so that
the regularity problem for $\CL \bu =0$ in $D_{\rho r}$ for $2<\rho<3$
with boundary data in $L^q_1$ is uniquely solvable \cite{G}.
It follows that
$$
\int_{\Delta_{2r}}
|(\nabla \bu)^*_{D_{2r}}|^q \, d\sigma
\le C\, \int_{\partial D_{\rho r}} |\nabla_t \bu|^q d\sigma
\le C\, \int_{\OO\cap \partial D_{\rho r}}
|\nabla \bu|^q \, d\sigma,
\tag 2.12
$$
where we have used the assumption that $\bu=0$ on $\Delta_{6r}$.
Integrating both sides of (2.12) in  $\rho \in (2, 3)$ yields that
$$
\frac{1}{|\Delta_{2r}|}\int_{\Delta_{2r}}
|(\nabla \bu)^*_{D_{2r}}|^q \, d\sigma
\le \frac{C}{|D_{3r}|}
\int_{D_{3r}} |\nabla \bu|^q dx.
\tag 2.13
$$
By Caccioppoli's
inequality as well as its well known consequence, the higher integrability
for $\nabla \bu$ \cite{Gi}, we obtain
$$
\aligned
&\left(\frac{1}{|D_{3r}|}
\int_{D_{3r}} |\nabla \bu|^q dx\right)^{1/q}
\le C\,
\left(\frac{1}{|D_{4r}|}
\int_{D_{4r}} |\nabla \bu|^2 dx\right)^{1/2}\\
&\le \frac{C}{r}
\left(\frac{1}{|D_{5r}|}
\int_{D_{5r}}
|\bu|^2 dx\right)^{1/2}
\le \frac{C}{r}
\left(\frac{1}{|\Delta_{5r}|}
\int_{\Delta_{5r}}
|(\bu)^*|^2 d\sigma\right)^{1/2}.
\endaligned
\tag 2.14
$$
Note that $q>2$ implies $p>2(n-1)/(n-3)$.
In view of (2.11), (2.13) and (2.14), we have proved that
for some $p>2(n-1)/(n-3)$,
$$
\aligned
\left(\frac{1}{|\Delta_r|}
\int_{\Delta_r}
|\Cal{M}_1(\bu)|^p d\sigma\right)^{1/p}
&
\le \sup_K |\bu|
+C\, \left(\frac{1}{|\Delta_{5r}|}
\int_{\Delta_{5r}} |(\bu)^*|^2 \, d\sigma\right)^{1/2}\\
&
\le 
C\, \left(\frac{1}{|\Delta_{5r}|}
\int_{\Delta_{5r}} |(\bu)^*|^2 \, d\sigma\right)^{1/2}.
\endaligned
\tag 2.15
$$
This, together with (2.8), gives the desired estimate (2.4).
\enddemo

To prove Theorem 1.3, it will be convenient to work with ``cubes'' on
$\partial\OO$. Let $D=\big\{ (x^\prime, x_n)\in \br^n:\  
x_n>\psi (x^\prime)\big\}$ and
$\partial D =\big\{ (x^\prime, \psi(x^\prime)):\ x^\prime\in\br^{n-1}
\big\}$. Define the map $\Phi: \partial D\to \br^{n-1}$
by $\Phi(x^\prime, \psi(x^\prime)) =x^\prime$.
We say that $Q\subset \partial D$ is a {\it cube} in $\partial D$, if
$\Phi(Q)$ is a cube in $\br^{n-1}$ with sides parallel to
the coordinate planes. A cube Q in $\partial D$ is said to be 
a dyadic subcube of $Q^\prime$ if $\Phi(Q)$ is a dyadic subcube
of $\Phi(Q^\prime)$ in $\br^{n-1}$, that is, $\Phi(Q)$
is one of the cubes obtained by bisecting the sides of $\Phi(Q^\prime)$
a finite number of times. Also, if $\Phi(Q)$ is a cube in $\br^{n-1}$,
$\rho \Phi(Q)$ denote the cube which has the same center, but
$\rho$ times the side length of $\Phi(Q)$. For a cube $Q$ on $\partial D$,
we denote $\Phi^{-1}\big[ \rho\Phi(Q)\big]$ by $\rho Q$.
As an example, in (2.2), $\Delta_{\rho r} =\rho\Delta_r$.

For cube $Q$ on $\partial D$ and a function $f$ defined 
on $Q$, we define a localized
Hardy-Littlewood maximal function $M_Q$ by
$$
M_Q (f) (P)
=\sup\Sb Q^\prime\owns P\\
Q^\prime\subset Q
\endSb
\frac{1}{|Q^\prime|}
\int_{Q^\prime}
|f|\,  d\sigma.
\tag 2.16
$$
In the next lemma, we will estimate $|E(\lambda)|$,
where
$$
E(\lambda)
=\big\{ P\in \Delta_r:\ M_{\Delta_{2r}} (|(\bu)^*|^2)(P)>\lambda\big\},\ \ \ 
0<r<c\, r_0.
\tag 2.17
$$

\proclaim{Lemma 2.18}
Suppose that $2<p<\frac{2(n-1)}{n-3} +\e$,
 where $\e>0$ is the same as in Lemma 2.3.
There exist positive constants $A$, $\delta$, $\gamma$ and $C_0$ depending
only on $n$, $m$, $\mu$, $p$ and $\OO$ such that if
$\CL\bu=\bo$ in $\OO$, $(\bu)^*\in L^2(\partial\OO)$ and
$\bu=\bbf\in L^2(\partial\OO)$ on $\partial\OO$, then
$$
|E(A\lambda)|
\le \delta |E(\lambda)|
+|\big\{
P\in \Delta_r:\ M_{\Delta_{2r}} (|\bbf|^2)(P)>\gamma \lambda\big\}|
\tag 2.19
$$
for all $\lambda\ge \lambda_0$, where
$$
\lambda_0=\frac{C_0}{|\Delta_{2r}|}\int_{\Delta_{2r}}
|(\bu)^*|^2 \, d\sigma.
\tag 2.20
$$
Most importantly, the constants $A$, $\delta$ satisfy 
the condition $\delta A^{p/2}<1$.
\endproclaim

\demo{Proof}
We use a real variable argument which is
  motivated by the method of approximation
 in \cite{CP}. We begin by fixing $p$ so that $2<p<2(n-1)/(n-3)+\e$.
Let $\delta\in (0,1)$ be a small constant to be determined.
By the $L^2$ boundedness of the Hardy-Littlewood maximal function,
$$
|E(\lambda)|
\le \frac{C(n, \|\psi\|_\infty)}{\lambda}
\int_{\Delta_{2r}} |(\bu)^*|^2 d\sigma.
\tag 2.21
$$
Thus, if $\lambda\ge\lambda_0$ where $\lambda_0$ is
 given in (2.20) with a large $C_0$,
we have $|E(\lambda)|<\delta |\Delta_r|$. Let $A=1/(2 \delta^{2/p})$.
Note that $\delta A^{p/2}=1/2^{p/2}<1/2$.

Since $E(\lambda)$ is  open relative to
$\Delta_r$, there exists a collection of disjoint dyadic subcubes
$\{ Q_k\}$ of $\Delta_r$ such that
$E(\lambda)=\bigcup_k Q_k$. Clearly we may 
assume that each $Q_k$ is maximal in the sense
that $\widetilde{Q}_k$ is not contained in $E(\lambda)$,
where $\widetilde{Q}_k$ denotes the dyadic ``parent'' of $Q_k$,
i.e., $Q_k$ is one of $2^{n-1}$ cubes obtained by
subdividing $\widetilde{Q}_k$ once.
Since $|E(\lambda)|<\delta |\Delta_r|$,
 we may also assume that
$32{Q}_k\subset \Delta_{2r}$ by taking $\delta $
sufficient small.

We claim that it is possible to choose positive constants
$\delta$, $\gamma$ and $C_0$ so that
if $\big\{ P\in {Q}_k: \ M_{\Delta_{2r}} (|\bbf|^2)(P)
\le \gamma \lambda\big\} \neq \emptyset$, then 
$$
|E(A \lambda)\cap Q_k|\le \delta |Q_k|.
\tag 2.22
$$
Estimate (2.19) follows from (2.22) by summation.

It remains to prove the claim. 
Suppose that $\big\{ P\in {Q}_k: \
M_{\Delta_{2r}} (|\bbf|^2)(P)\le \gamma \lambda\big\}
\neq \emptyset$. Since $Q_k$ is maximal,
 by a simple geometric observation,
we obtain
$$
M_{\Delta_{2r}} (|(\bu)^*|^2)(P)
\le \max\big(
M_{2 Q_k} (|(\bu)^*|^2)(P), C_1\, \lambda\big)
\tag 2.23
$$
for any $P\in Q_k$, where $C_1$ depends only on $n$ and 
$\| \nabla\psi\|_\infty$. We may assume that $A\ge C_1$
by  making $\delta$ small.
It follows that
$$
|Q_k\cap E(A\lambda)|
\le |\left\{ P\in Q_k:\ 
M_{2{Q}_k} (|(\bu)^*|^2)(P)>A\lambda\right\}|.
\tag 2.24
$$
Also note that, if $\widetilde{Q}_k\subset Q\subset \Delta_{2r}$, then
$$
\frac{1}{|Q|}
\int_Q |\bbf|^2 \, d\sigma \le \gamma \lambda\ \text{ and }\ 
\frac{1}{|Q|}
\int_Q |(\bu)^*|^2 \, d\sigma \le \lambda.
\tag 2.25
$$

Let $\bold{v}=\bold{v}_k$ be the unique solution of the 
$L^2$ Dirichlet problem
on $\OO$  with boundary data $\bbf\chi_{16{Q}_k}$. It follows from
(2.24) and Lemma 2.3 that
$$
\aligned
&|Q_k\cap E(A\lambda)|
\le |\left\{ P\in Q_k: M_{2 {Q}_k}
(|(\bu-\bold{v})^*|^2)(P)>\frac{A\lambda}{4}\right\}| \\
&\ \ \ \ \ \ \ \ \ \ \ \ \ \ \ \ \ \ 
+
|\left\{ P\in Q_k: M_{2 {Q}_k}
(|(\bold{v})^*|^2)(P)>\frac{A\lambda}{4}\right\}|\\
&
\le \frac{C}{(A\lambda)^{\bar{p}/2}}
\int_{2 {Q}_k} |(\bu-\bold{v})^*|^{\bar{p}} \, d\sigma
+\frac{C}{A\lambda}
\int_{2 {Q}_k}
|(\bold{v})^*|^2 \, d\sigma\\
&
\le |Q_k|
\left\{
C\, 
\left(\frac{1}{A\lambda |10 {Q}_k|}
\int_{10 {Q}_k}
|(\bu-\bold{v})^*|^2 \, d\sigma\right)^{\bar{p}/2}
+\frac{C}{A\lambda |16 {Q}_k|}
\int_{16 {Q}_k}|\bbf|^2 \, d\sigma\right\},
\endaligned
$$
where $p<\bar{p}<2(n-1)/(n-3)+\e$ and
 we also used the $L^2$ estimate for $\bold{v}$ in the last inequality.
This, together with (2.25), the $L^2$ estimate and 
the choice $A=1/(2\delta)^{2/p}$,  gives
$$
|Q_k\cap E(A\lambda)|
\le |Q_k|
\left\{ \frac{C_2}{A^{\bar{p}/2}}
+\frac{C_2\,\gamma}{A}\right\}
= \delta |Q_k|
\left\{ C_2 2^{\bar{p}/2} \delta^{\frac{\bar{p}}{p}-1}
+2 C_2 \gamma \delta^{\frac{2}{p}-1}\right\},
\tag 2.26
$$
where $C_2$ depends only on $n$, $m$, $\mu$, $\bar{p}$ and $\OO$.

Finally, since $\bar{p}>p$, we may choose $\delta>0$ so small that
$C_2 2^{\bar{p}/2} \delta^{\frac{\bar{p}}{p}-1}<1/2$. 
With $\delta$ fixed, we choose $\gamma>0$ so small that
$2 C_2 \gamma \delta^{\frac{2}{p}-1}<1/2$.
We then obtain $|Q_k\cap E(A\lambda)|<\delta |Q_k|$.
The proof is complete.
\enddemo

Estimate (2.19) is a good-$\lambda$ type inequality.
Its relation with the $L^p$ norm estimate is well known.

\proclaim{Lemma 2.27}
Suppose that $2<p<\frac{2(n-1)}{n-3}+\e$ where $\e>0$
is the same as in Lemma 2.3.
If $\bbf\in L^p(\partial\OO)$ and $\bu$ is the unique
solution of the $L^2$ Dirichlet problem with boundary data $\bbf$, then
$(\bu)^*\in L^p(\Delta_r)$ and 
$$
\aligned
&\left(\frac{1}{|\Delta_r|}
\int_{\Delta_r}
|(\bu)^*|^p \, d\sigma\right)^{1/p}\\
&\le C\, 
\left(\frac{1}{|\Delta_{2r}|}
\int_{\Delta_{2r}}
|(\bu)^*|^2 \, d\sigma\right)^{1/2}
+C\, \left(\frac{1}{|\Delta_{2r}|}
\int_{\Delta_{2r}}
|\bbf|^p \, d\sigma\right)^{1/p}.
\endaligned
\tag 2.28
$$
\endproclaim

\demo{Proof}
We begin by multiplying both sides of (2.19) by $\lambda^{\frac{p}{2}-1}$
and integrating the resulting inequality in $\lambda$
on the interval  $(\lambda_0,\Lambda)$. This gives
$$
\frac{1}{A^{p/2}}
\int_{A\lambda_0}^{A\Lambda}
\lambda^{\frac{p}{2}-1} |E(\lambda)| \, d\lambda
\le \delta \int_{\lambda_0}^\Lambda\lambda^{\frac{p}{2}-1}
|E(\lambda)|\, d\lambda
+C\int_{\Delta_{2r}} |\bbf|^p d\sigma.
\tag 2.29
$$
It follows that
$$
\aligned
\left(\frac{1}{A^{p/2}}-\delta\right)
\int_0^\Lambda
\lambda^{\frac{p}{2}-1} |E(\lambda)| \, d\lambda
&\le C\int_0^{A\lambda_0}
\lambda^{\frac{p}{2}-1} |E(\lambda)| \, d\lambda
+C\, \int_{\Delta_{2r}}|\bbf|^p \, d\sigma\\
&\le C\lambda_0^{\frac{p}{2}}
|\Delta_{2r}|
+C\, \int_{\Delta_{2r}}|\bbf|^p \, d\sigma.
\endaligned
$$
Since $\delta A^{p/2}<1$, estimate (2.28)
follows easily from the above inequality and (2.20)
by letting $\Lambda\to\infty$.
\enddemo

We are now in a position to give the proof of Theorem 1.3.

\demo{\bf Proof of Theorem 1.3}
The case $2-\e<p\le 2$ is already known \cite{G}.
For $p>2$, 
since $L^p(\partial\OO)\subset L^2(\partial\OO)$ for $p>2$, the 
uniqueness follows from the uniqueness for the
Dirichlet problem with $L^2$ data. The existence as well as the estimate
$\|(\bu)^*\|_p\le C\, \| \bbf\|_p$ follows easily from Lemma 2.27
by covering $\partial\OO$ with a finite number of
coordinate patches.
\enddemo

\remark{\bf Remark 2.30}
Theorem 1.3 also holds for the exterior domain $\OO_-
=\br^n\setminus \overline{\OO}$. That is, given any
$\bbf\in L^p(\partial\OO)$ with $2\le p<2(n-1)/(n-3)+\e$,
there exists a unique $\bu$ on $\OO_-$ satisfying
$ \CL \bu =\bo$ in $\OO_-$, $\bu=\bbf$ on $\partial\OO$,
$(\bu)^*_e\in L^p(\partial\OO)$, and 
$|\bu(x)|=O(|x|^{2-n})$ as $|x|\to\infty$.
Moreover, the solution $\bu$ satisfies
$\| (\bu)^*_e\|_p\le C\, \| \bbf\|_p$, where
 $(\bu)^*_e$ denotes the nontangential maximal
function of $\bu$ with respect to $\OO_-$.
\endremark

\remark{\bf Remark 2.31}
In the case $n=3$,  the techniques in  this section
yield the $L^p$ solvability of the Dirichlet
problem for $2<p<\infty$.
\endremark

We end this section by formulating a theorem, which is
of independent interest, on the $L^p$ boundedness
of operators on $\br^n$. Its proof, which is given in \cite{S5},
may be carried out by slightly modifying the techniques
used in the proof of Lemma 2.18.
We point out that any operator with a standard Calder\'on-Zygmund
kernel satisfies condition
(2.33) for all $p>2$.

\proclaim{Theorem 2.32}
Let $T$ be a bounded sublinear operator on $L^2(\br^n)$. 
Suppose that for some $p>2$,  $T$ satisfies the following
$L^p$ localization property:
$$
\aligned
&\left\{
\frac{1}{|Q|}
\int_{Q} | Tf|^p\, dx\right\}^{1/p}
\\
&\le C\,
\left\{
\bigg(\frac{1}{|2Q|}\int_{2Q}
|Tf|^2\, dx\bigg)^{1/2}
+\sup_{Q^\prime\supset Q}
\bigg(\frac{1}{|Q^\prime|}
\int_{Q^\prime}
|f|^2\, dx\bigg)^{1/2}\right\},
\endaligned
\tag 2.33
$$
for any cube $Q\subset\br^n$ and any
$L^2$ function $f$ with supp$(f)\subset \br^n\setminus 3Q$.
Then $T$ is bounded on $L^q(\br^n)$ for
any $2<q<p$.
\endproclaim

\medskip

\centerline{\bf 3. The Polyharmonic Equation}

In this section we give the proof of Theorem 1.10, using the
same line of argument as in the last section.
We will need the following estimates for suitable solutions of
$\Delta^\ell u=0$ on $\OO$:
$$
\align
\| (\nabla^{\ell-1} u)^*\|_2 &\le C\, \| \nabla ^{\ell -1}u\|_2,
\tag 3.1\\
\| (\nabla^\ell u)^*\|_q
&\le C\, \| \nabla_t \nabla^{\ell-1} u\|_q,
\tag 3.2
\endalign
$$
for some $q>2$ depending only on $n$, $\ell$, and the
Lipschitz character of $\OO$. Both estimates
were established in \cite{V3} for any integer $\ell\ge 2$.
In the case of the biharmonic equation ($\ell=2$),
 estimate (3.1) was obtained
earlier in \cite{DKV1}.
We also need a Caccioppoli's inequality,
$$
\int_{D_{sr}} |\nabla^\ell u|^2\, dx
\le \frac{C}{(\rho -s)^2 r^2}
\int_{D_{\rho r}} |\nabla^{\ell-1} u|^2 \, dx,
\tag 3.3
$$
where $0<s<\rho<1$, $u$ satisfies $\Delta^\ell u =0$ in $\OO$,
$(\nabla ^\ell u)^*\in L^2 (\Delta_{2r})$, and
$D^\alpha u=0$ on $\Delta_{3r}$ for
all $0\le |\alpha|\le \ell-1$.
Inequality (3.3) may be proved by using integration by parts
and Poincar\'e's inequalities for functions which
vanish on part of the boundary. From 
(3.3) one may deduce the following boundary reverse H\"older inequality
for some $q>2$
by a standard argument,
$$
\left(\frac{1}{|D_r|} \int_{D_r}
|\nabla^\ell u|^q dx\right)^{1/q}
\le C\left(\frac{1}{|D_{2r}|}
\int_{D_{2r}} |\nabla^\ell u|^2 \, dx\right)^{1/2}.
\tag 3.4
$$ 

\proclaim{Lemma 3.5}
There exists $\e>0$ depending only on $n$, $\ell$ and $\OO$ such that if
$\Delta^\ell u=0$ in $\OO$, $(\nabla ^{\ell-1} u)^*\in L^2(\partial\OO)$
and $D^\alpha u =0$ for all $0\le |\alpha|\le
\ell-1$ on the surface cube $\Delta_{8r}$, then
$(\nabla^{\ell-1} u)^*\in L^p(\Delta_r)$ and
$$
\left(\frac{1}{|\Delta_r|}
\int_{\Delta_r}|(\nabla^{\ell -1} u)^*|^p \, d\sigma\right)^{1/p}
\le C\,
\left(\frac{1}{|\Delta_{5r}|}
\int_{\Delta_{5r}} |(\nabla^{\ell -1} u)^*|^2 \, d\sigma\right)^{1/2},
\tag 3.6
$$
where $2<p<\frac{2(n-1)}{n-3} +\e$.
\endproclaim

\demo{Proof}
Using the well known interior estimates, one may prove that for $P\in\Delta_r$,
$$
\aligned
&(\nabla^{\ell -1} u)^*(P)\\
&\le C\left(\frac{1}{|\Delta_{5r}|}
\int_{\Delta_{5r}}
|(\nabla^{\ell-1} u)^*|^2 \, d\sigma\right)^{1/2}
+C\, \int_{\Delta_{2r}}
\frac{(\nabla^\ell u)^*_{D_{2r}}(Q)}{
|P-Q|^{n-2}}\, d\sigma (Q).
\endaligned
$$
With estimates (3.2), (3.4) and (3.3) at our
disposal, the rest of the proof is similar to that of
Lemma 2.1. We omit the details.
\enddemo

As in (2.17), we let
$$
F(\lambda)
=\left\{ P\in \Delta_r:\  M_{\Delta_{2r}}
(|(\nabla^{\ell-1} u)^*|^2)(P)>\lambda\right\} \ \ \text{ for } 0<r<c\, r_0.
\tag 3.7
$$

\proclaim{Lemma 3.8}
Suppose that $2<p<\frac{2(n-1)}{n-3}+\e$ where $\e>0$ is the same as in
Lemma 3.5. There exist positive constants $A$, $\delta$, $\gamma$ and
$C_0$ depending only on $n$, $\ell$, $p$ and $\OO$, such that
if $f=\{ f_\alpha: \, 0\le |\alpha|\le 
\ell-2\} \in WA^2_{\ell-1}
(\partial\OO)$, $g\in L^2(\partial\OO)$ and $u$
is the unique solution of the $L^2$ Dirichlet problem (1.9)
with boundary data $f$, $g$, then
$$
|F(A\lambda)|
\le \delta |F(\lambda)|
+\big|\bigg\{ P\in\Delta_r:\
M_{\Delta_{2r}} \big\{|g|^2+\sum_{|\alpha|=\ell-2}
|\nabla_t f_\alpha|^2\big\}(P)>\gamma \lambda\bigg\}\big|
\tag 3.9
$$
for all $\lambda\ge \lambda_0$, where
$$
\lambda_0=\frac{C_0}{|\Delta_{2r}|}
\int_{\Delta_{2r}}
|(\nabla^{\ell -1} u)^*|^2 \, d\sigma.
\tag 3.10
$$
Most importantly, the constants $A$, $\delta$ satisfy the condition
$\delta A^{p/2}<1$.
\endproclaim

\demo{Proof}
The proof is similar to that of Lemma 2.18.
However, in the place of $\bold{v}$, we have to be
more careful with the choice of
polyharmonic function $v=v_k$ in $\OO$ for each $Q_k$,
since $WA^2_{\ell-1}(\partial\OO)$ is a Sobolev space.
Let 
$\varphi=\varphi_k$ is a smooth cut-off function on $\br^n$
 such that
$\varphi=1$ on $16{Q}_k$, $\varphi =0$ on 
$\partial\OO\setminus 17{Q}_k$, $|D^\alpha \varphi|
\le C/r^{|\alpha|}$ for $0\le |\alpha|\le \ell-1$,
where $r=r_k$ is the diameter of $Q_k$.
Let $h$ be a polynomial of degree $\ell -2$ 
to be determined. We choose $v=v_k$ to be the solution 
of $L^2$ Dirichlet problem (1.9) with
boundary data 
$$
D^\alpha v
=D^\alpha \big((u-h)\varphi\big)
=\sum_{0\le \beta\le \alpha}
\frac{\alpha!}{\beta! (\alpha-\beta)!}
(f_\beta -D^\beta h) D^{\alpha-\beta} \varphi 
\tag 3.11
$$ for $0\le |\alpha|
\le \ell -2$ and  $\frac{\partial^{\ell -1} v}
{\partial N^{\ell -1}} =g\varphi$ on $\partial\OO$.

Write $u=v+w +h$. Then $\nabla^{\ell -1} u
=\nabla^{\ell-1} v +\nabla^{\ell -1} w$, since $h$ is a polynomial
of degree $\ell -2$. It is easy to see that
 $w$ is the solution
of the $L^2$ Dirichlet problem (1.9) with boundary values
vanishing on $16{Q}_k$. In fact,
$D^\alpha w =D^\alpha \big((u-h)(1-\varphi)\big)$
for $0\le |\alpha|\le \ell-2$
and $\frac{\partial^{\ell-1} w}{\partial N^{\ell-1}}
=g(1-\varphi)$ on $\partial\OO$.
Thus we may apply Lemma 3.5 to $w$.

To control the $L^2$ norm of $(\nabla^{\ell-1} v)^*$, we 
 note that for $|\alpha|=\ell-2$
$$
\aligned
&\| \nabla_t \nabla^\alpha \big((u-h)\varphi)\|_2\\
&\le
C\, \sum_{0\le \beta\le \alpha}
\bigg\{ \| \big(\nabla_t (f_\beta -D^\beta h)\big) 
D^{\alpha-\beta} \varphi\|_2
 +
\| (f_\beta-D^\beta h)\nabla_t D^{\alpha-\beta}\varphi\|_2\bigg\}\\
&\le
C\, \sum_{0\le \beta\le \alpha}
r^{|\beta|-\ell+2}\left\{
\|\nabla_t (f_\beta -D^\beta h )\|_{L^2(17{Q}_k)}
+r^{-1}
\| f_\beta -D^\beta h\|_{L^2(17{Q}_k)}\right\}.
\endaligned
\tag 3.12
$$

Finally, we let $h(x)=\sum_{|\alpha|\le \ell -2} \frac{B_\alpha}
{\alpha!}\, x^\alpha$, where the constants $B_\alpha$ are defined
inductively by
$$
\aligned
&
B_\alpha =\frac{1}{|17{Q}_k|}
\int_{17{Q}_k} f_\alpha (P) \, d\sigma(P)\ \text{ for }
|\alpha|=\ell-2,\\
&
B_\alpha
=\frac{1}{|17{Q}_k|}\int_{17{Q}_k}
\bigg\{
f_\alpha (P)-\sum\Sb |\alpha|\le \ell -2\\
\alpha>\beta
\endSb \frac{B_\alpha}{(\alpha -\beta)!}\, P^{\alpha -\beta}\bigg\}
\, d\sigma(P)
\endaligned
\tag 3.13
$$
for $0\le |\alpha|<\ell-2$.
We remark that $B_\alpha$ is defined in such a way that
$$
\int_{17{Q}_k} \big\{ f_\beta -D^\beta h\big\}\, d\sigma
=0\ \ \text{for all }0\le |\beta|\le \ell-2.
\tag 3.14
$$
From this, by using
 Poincar\'e's inequality repeatedly, we obtain
$$
\aligned
\| f_\beta -D^\beta h\|_{L^2(17{Q}_k)}
&\le C\, r \, \sum_{|\alpha|=|\beta| +1}
\| f_\alpha -D^\alpha h\|_{L^2(17{Q}_k)}\\
&\le
C\, r^{\ell -|\beta|-2}
\sum_{|\alpha|=\ell -2}
\| f_\alpha -D^\alpha h\|_{L^2(17{Q}_k)}\\
&\le C\, r^{\ell-|\beta|-1}
\sum_{|\alpha|=\ell-2}
\|\nabla_t f_\alpha\|_{L^2(17{Q}_k)}.
\endaligned
\tag 3.15
$$
By combining estimates (3.12) with (3.15), we get
$$
\| (\nabla^{\ell -1} v)^*\|_2
\le C\, \bigg\{
\| g\|_{L^2(17{Q}_k)}
+\sum_{|\alpha|=\ell -2}
\|\nabla_t f_\alpha \|_{L^2(17{Q}_k)}
\bigg\}.
\tag 3.16
$$
With these observations, the argument in Lemma 2.18 goes through
with minor changes. We leave the details to the reader.
\enddemo

The following lemma follows from Lemma 3.4 by integration.

\proclaim{\bf Lemma 3.17}
Suppose that $2<p<\frac{2(n-1)}{n-3} +\e$, where $\e>0$ is the same as
in Lemma 3.5. If $f\in WA^p_{\ell-1}(\partial\OO)$, $g\in L^p(\partial\OO)$,
and $u$ is the polyharmonic function on $\OO$ 
satisfying $(\nabla^{\ell-1} u)^*\in L^2(\partial\OO)$
and $D^\alpha u=f_\alpha$ for $0\le |\alpha|\le \ell-2$, 
$\frac{\partial^{\ell-1} u}{\partial N^{\ell-1}} =g$ on
$\partial\OO$, then $(\nabla^{\ell-1} u)^*\in L^p(\Delta_r)$ and 
$$
\aligned
&\left\{\frac{1}{|\Delta_r|}
\int_{\Delta_r}
|(\nabla^{\ell-1} u)^*|^p \, d\sigma\right\}^{1/p}
\le C\, 
\left\{\frac{1}{|\Delta_{2r}|}
\int_{\Delta_{2r}}
|(\nabla^{\ell-1} u)^*|^2 \, d\sigma\right\}^{1/2}\\
&\ \ \ \ \ \ \ \ \ \ \ \ +C\, \bigg\{\frac{1}{|\Delta_{2r}|}
\int_{\Delta_{2r}}
\big(|g|^p+\sum_{|\alpha|=\ell -2} |\nabla_t f_\alpha|^p
\big) \, d\sigma\bigg\}^{1/p}.
\endaligned
\tag 3.18
$$
\endproclaim

As before, Theorem 1.10 is a consequence of the $L^2$ solvability
and Lemma 3.17.

\medskip

\centerline{\bf 4. The Regularity Problem for Elliptic Systems}

This section is devoted to the proof of Theorem 1.4. 
We follow an approach found in \cite {V1}, where it was used by
Verchota to establish the solvability of the $L^p$ regularity
problem  for Laplace's equation.
 The basic idea is to prove
$$
\| \frac{\partial \bu}{\partial\nu}\|_p \le C\, \| \nabla_t \bu\|_p,
\tag 4.1
$$
for starshaped Lipschitz domains,
using a duality argument which involves conjugate functions
and area integral estimates.
With (4.1) one may show that the single layer potential,
which maps $L^p(\partial\OO)$ to $L^p_1(\partial\OO)$, is
invertible.
We should mention that in (4.1), $\frac{\partial \bu}{\partial\nu}$
denotes the conormal derivative defined by
$\big(\frac{\partial \bu}{\partial\nu}\big)^r
=a_{ij}^{rs} \frac{\partial u^s}{\partial x_j} N_i$.

\proclaim{Lemma 4.2}
Let $\OO$ be a starshaped Lipschitz domain in $\br^n$, $n\ge 3$.
 Suppose that for some $p>1$, the $L^{p}$ Dirichlet problem (1.2) is 
uniquely solvable for any $\bbf\in L^p(\partial\OO)$.
If $\CL \bu =\bo$ in $\OO$, $(\nabla \bu)^*
\in L^{p^\prime}(\partial\OO)$ and $\nabla \bu$ has nontangential
limit on $\partial\OO$, then
$\| \frac{\partial \bu}{\partial \nu}\|_{p^\prime}
\le C\, \| \nabla_t \bu\|_{p^\prime}$
where $p^\prime=p/(p-1)$.
\endproclaim

\demo{Proof}
We may assume that $\OO$ is starshaped with respect to the
origin.
Let $\bg$ be a Lipschitz continuous function on $\partial\OO$
and $\bv$ be the solution of $\CL \bold{v}=0$
in $\OO$ satisfying $(\nabla \bv)^*\in L^2(\partial\OO)$ and $\bv=\bg$ on $
\partial\OO$.
Define
$$
\bold{H}(x)=\int_0^1  \bw (rx)
\, \frac{dr}{r}\ \ \text{ for } x\in \OO,
\tag 4.3
$$
where $\bw(x)=\bv(x)-\bv(0)$.
It is easy to verify that $\CL \bold{H} =\bo$ and
$\frac{\partial H^r}{\partial x_i} \, x_i =w^r$.
As in the case of harmonic functions, we also have
$$
\| (\nabla \bold{H})^*\|_p
 \le C\, \|(\bw)^*\|_p
\le C\, \| \bg\|_p.
\tag 4.4
$$
The first inequality in (4.4), which holds for any $p>0$,
follows from 
 the area integral estimates
for elliptic systems \cite{DKPV}, and  the second 
 by the solvability of the
$L^p$ Dirichlet problem. We refer the reader to \cite{V2, p870-871}
for details in the case of Laplace's equation.

We now use integration by parts to obtain
$$
\aligned
\int_{\partial\OO} \frac{\partial \bu}{\partial\nu} \cdot \bg\, d\sigma
&=\int_{\partial\OO} \bu \cdot \frac{\partial \bv}{\partial \nu}\, d\sigma
=\int_{\partial\OO} \bu \cdot \frac{\partial \bw}{\partial \nu}\, d\sigma\\
&
=(2-n)\int_{\partial\OO}
\bu \cdot \frac{\partial \bold{H}}{\partial \nu}\, d\sigma
-\int_{\partial\OO}
\frac{\partial u^r}{\partial t_{i\ell}}
a_{ij}^{rs} \frac{\partial H^s}{\partial x_j} x_\ell \, d\sigma,
\endaligned
\tag 4.5
$$
where $\frac{\partial }{\partial t_{i\ell}}
=n_i\frac{\partial}{\partial x_\ell}
-n_\ell \frac{\partial}{\partial x_i}$ is a tangential
derivative.
We remark that the integration by parts used in (4.5) can be 
justified by an approximation argument, using estimate (4.4) and
the assumption that $\nabla \bu$ has nontangential limits on 
$\partial\OO$ and $(\nabla \bu)^*
\in L^{p^\prime}$.
It follows from (4.5) and H\"older's inequality  that
$$
\left|
\int_{\partial\OO}
\frac{\partial\bu}{\partial\nu} \cdot \bg \, d\sigma\right|
\le C\, \| \nabla _t \bu\|_{p^\prime}
\|\bg\|_p,
\tag 4.6
$$
where (4.4) and Poincar\'e's inequality are also used.
The desired estimate then follows by duality.
\enddemo

We now give the proof of Theorem 1.4.

\demo{\bf Proof of Theorem 1.4}
Let $\Gamma(x)$ denote the matrix fundamental solution for the operator
$\CL$ on $\br^n$ with pole at the origin. Consider the
single layer potential
$$
\Cal{S}(\bg)(x)=\int_{\partial\OO} \Gamma(x-y) \bg(y)\, d\sigma.
\tag 4.7
$$
For the existence as well as the estimate $\|(\nabla \bu)^*\|_p
\le C\| \nabla_t \bu\|_p$, it suffices to show that there exists $\e>0$ such 
that
$\Cal{S}: L^p(\partial\OO)\to L^p_1(\partial\OO)$
is invertible for $\frac{2(n-1)}{n+1} -\e <p\le 2$.
To this end, we fix $P_0\in \partial\OO$ and assume that
$$
B(P_0,r_0)\cap \OO
=B(P_0, r_0)\cap \big\{ (x^\prime, x_n)\in \br^n:\ \ 
x_n>\psi(x^\prime)\big\}
\tag 4.8
$$
where $\psi:\br^{n-1}\to\br$ is Lipschitz continuous.
We may assume that $P_0=0$.
Let
$$
\OO_r =
\big\{ (x^\prime, x_n)\in \br^n:\ \ 
|x_1|<r, \dots, |x_{n-1}|<r, \, \psi(x^\prime)<x_n
<B\, r\big\}
\tag 4.9
$$
where the constant $B=B(n,\|\nabla\psi\|_\infty)>0$ is chosen
so that $\OO_r$ is a starshaped Lipschitz domain for any $r>0$.
Note that, by Theorem 1.3, there exists $\delta>0$ depending only on
$n$, $m$, $\mu$ and $\| \nabla \psi\|_\infty$, such that
the $L^q$ Dirichlet problem for the operator $\CL$ on
$\OO_r$ is uniquely solvable for any $r>0$
and $2\le q <q_0=\frac{2(n-1)}{n-3} +\delta$.

Let $\bg\in L^p(\partial\OO)$ and
 $\bv =\Cal{S}(\bg)$ in $\br^n$. Suppose $2\le p^\prime<q_0$.
We may apply Lemma 4.2 to $\bv$
on $\OO_r$ for $0<r<c\, r_0$. This gives
$$
\int_{\Delta_r} \big|\frac{\partial\bv}{\partial \nu}\big|^p d\sigma
\le C\, \int_{\partial \OO_{sr}}
|\nabla_t \bv|^p \, d\sigma
\le C\, \int_{\partial\OO}
|\nabla_t \bv|^p\, d\sigma
+C\, \int_{\OO\cap \OO_{sr}} |\nabla \bv|^p\, d\sigma,
\tag 4.10
$$
where $s\in (1,2)$. By integrating (4.10) in $s$ and 
covering $\partial\OO$ with a finite number of
coordinate patches, we obtain
$$
\aligned
\int_{\partial\OO}
\big|\frac{\partial \bv}{\partial\nu}\big|^p\, d\sigma
&\le C\,
\int_{\partial \OO}
|\nabla_t \bv|^p \, d\sigma
+C\, \int_{\widetilde{\OO}} |\nabla \bv|^p \, dx\\
& 
\le C\,
\int_{\partial \OO}
|\nabla_t \bv|^p \, d\sigma
+\gamma \int_{\partial\OO} |(\nabla \bv)^*|^p \, d\sigma
+C_\gamma \sup_K |\bv|^p,
\endaligned
\tag 4.11
$$
where $\widetilde{\OO}=\{ x\in \OO: \ \text{dist}(x,
\partial\OO)<r_0\}$, and $K$ is a compact set in $\OO$.
To estimate $\bv$ on $K$, we choose $\bar{q}<2$ so that
the $L^{\bar{q}}$ Dirichlet problem is uniquely solvable.
It follows that
$$
\sup_K |\bv|\le C\, \|(\bv)^*\|_{\bar{q}}
\le C\, \| \bv\|_{\bar{q}}
\le C\, \big\{ \|\nabla_t \bv\|_{\bar{p}}
+\| \bv\|_{\bar{p}}\big\},
\tag 4.12
$$
where $\frac{1}{\bar{q}} =\frac{1}{\bar{p}}-\frac{1}{n-1}$,
and we used the Sobolev imbedding in the last inequality.
Note that $\bar{q}<2$ implies $\bar{p}<2(n-1)/(n+1)$.
We may assume that $\bar{p}^\prime<2(n-1)/(n-3) +\delta$.
Thus, if   $\bar{p}\le p\le 2$, we obtain
$$
\int_{\partial\OO}
\big|\frac{\partial \bv}{\partial\nu}\big|^p\, d\sigma
\le C\,
\int_{\partial \OO}
|\nabla_t \bv|^p \, d\sigma
+C\, \int_{\partial\OO}
|\bv|^p\, d\sigma
+\gamma \, C\, \int_{\partial \OO} |\bg|^p d\sigma,
\tag 4.13
$$
where we have used the estimate $\|(\nabla \bv)^*\|_p
\le C\, \| \bg\|_p$ (see e.g. \cite{V1} or \cite{K1}).
By the same argument, estimate (4.13) also holds for
the exterior domain $\br^n\setminus \overline{\OO}$.
Thus,  by the jump 
relation $\bg =\frac{\partial \bv_+}{\partial\nu}
-\frac{\partial \bv_-}{\partial\nu}$ where
$\pm$ indicate the limits taken from $\OO$ and
$\br^n\setminus\overline{\OO}$
respectively,
we obtain
$$
\| \bg\|_p \le \|\frac{\partial \bv_+}{\partial\nu}\|_p
+\|\frac{\partial \bv_-}{\partial\nu}\|_p
\le C\,
\big\{
\| \nabla \bv_t\|_p +\| \bv\|_p +\gamma \| \bg\|_p\big\},
\tag 4.14
$$
for $\bar{p}\le p\le 2$. Choose $\gamma>0$ small so that
$\gamma C<1/2$. This gives
$$
\|\bg\|_p\le C\, \| \Cal{S}(\bg)\|_{L^p_1(\partial\OO)}.\tag 4.15
$$
Estimate (4.15) implies that
the operator $\Cal{S}:L^p(\partial\OO)\to L^p_1(\partial\OO)$ 
is one-to-one and has a closed range. Note that the range
is also dense, since $\Cal{S}: L^2(\partial\OO)\to L^2_1(\partial\OO)$
is invertible \cite{G}. We conclude that 
$\Cal{S}: L^p(\partial\OO)\to L^p_1(\partial\OO)$
is invertible for $\bar{p}\le p\le 2$.

Finally, to show the uniqueness, we suppose that $\CL{\bu}=\bo$
in $\OO$, $\bu=0$ on $\partial \OO$ and $(\nabla \bu)^*
\in L^p$ for some $p>2(n-1)/(n+1)-\e$.
It follows from (2.10) that $(\bu)^*\in L^q$, where $\frac{1}{q}
=\frac{1}{p}-\frac{1}{n-1}$.
Note that if $\e>0$ is small, $q>2-\bar{\e}$.
Thus the uniqueness in Theorem 1.4 follows from the uniqueness
of the $L^q$ Dirichlet problem for $q>2-\bar{\e}$ \cite{G}.
\enddemo


\Refs
\widestnumber\key{DKPV}

\ref\key CP
\by L.A.~Caffarelli and I.~Peral
\paper On $W^{1,p}$ estimates for elliptic equations
in divergence form
\jour Comm. Pure App. Math. \vol 51
\yr 1998
\pages 1-21
\endref

\ref \key  D1
\by B.~Dahlberg
\paper On estimates for harmonic measure
\jour Arch. Rat. Mech. Anal. \vol 65
\yr 1977 \pages 273-288
\endref


\ref \key D2
\by B.~Dahlberg
\paper On the Poisson integral for Lipschitz and $C^1$ domains
\jour Studia Math. \vol 66
\yr 1979 \pages 13-24
\endref

\ref \key DK1
\by B.~Dahlberg and C.~Kenig
\paper Hardy spaces and the Neumann problem in $L^p$
for Laplace's equation in Lipschitz domains
\jour Ann. of Math. \vol 125\yr 1987
\pages 437-466
\endref


\ref\key DK2
\by B.~Dahlberg and C.~Kenig
\paper $L^p$ estimates for the three-dimensional systems of
elastostatics on Lipschitz domains
\jour Lecture Notes in Pure and Appl. Math.
\vol 122 \yr 1990 \pages 621-634
\endref

\ref \key DKPV
\by B.~Dahlberg, C.~Kenig, J.~Pipher, and G.~Verchota
\paper Area integral estimates for higher order elliptic
equations and systems
\jour Ann. Inst. Fourier (Grenoble) \vol 47 \pages 1425-1461
\yr 1997
\endref

\ref \key DKV1
\by B.~Dahlberg, C.~Kenig, and G.~Verchota
\paper
The Dirichlet problem for the biharmonic equation
in a Lipschitz domain
\jour Ann. Inst. Fourier (Grenoble)
\vol 36 \yr 1986 \pages 109-135
\endref


 
\ref\key DKV2
\by B.~Dahlberg, C.~Kenig, and G.~Verchota
\paper 
Boundary value problems for the systems of elastostatics
in Lipschitz domains
\jour Duke Math. J.
\vol 57 \yr 1988 \pages 795-818
\endref

\ref \key F
\by E.~Fabes
\paper
Layer potential methods for boundary value problems on 
Lipschitz domains
\jour Lecture Notes in Math.
\vol 1344
\yr 1988
\pages 55-80
\endref

\ref\key FKV
\by E.~Fabes, C.~Kenig, and G.~Verchota
\paper
Boundary value problems for the Stokes system on Lipschitz
domains
\jour Duke Math. J. \vol 57 \yr 1988
\pages 769-793
\endref


\ref \key G
\by W.~Gao 
\paper Boundary value problems on Lipschitz domains for
general elliptic systems
\jour J. Funct. Anal.
\yr 1991\pages 377-399
\endref

\ref\key Gi
\by M.~Giaquinta
\book Multiple Integrals in the Calculus of Variations
and Nonlinear Elliptic Systems
\publ Annals of Math. Studies 105,
Princeton Univ. Press
\yr 1983
\endref

\ref \key JK
\by D.~Jerison and C.~Kenig
\paper The Neumann problem in Lipschitz domains
\jour Bull. Amer. Math. Soc.
\vol 4 \yr 1981 \pages 203-207
\endref

\ref \key K1
\by C.~Kenig
\paper Elliptic boundary value problems on Lipschitz domains
\jour Beijing Lectures in Harmonic Analysis,
Ann. of Math. Studies
\vol 112
\yr 1986
\pages 131-183
\endref



\ref\key K2
\by C.~Kenig
\book Harmonic Analysis Techniques for Second
Order Elliptic Boundary Value Problems
\bookinfo
CBMS Regional Conference Series in Math.\vol 83
\publ AMS, Providence, RI
\yr 1994
\endref

\ref\key MM
\by D.~Mitrea and M.~Mitrea,
\paper General second order, strongly elliptic
systems in low dimensional nonsmooth manifold
\jour
Contemporary Math. 
\vol 277 \yr 2001
\pages 61-86
\endref

\ref\key PV1
\by J.~Pipher and G.~Verchota
\paper The Dirichlet problem in $L^p$ for the
biharmonic equation on Lipschitz domains
\jour Amer. J. Math. \vol 114 \yr 1992 \pages 923-972
\endref


\ref\key  PV2
\by J.~Pipher and G.~Verchota
\paper A maximum principle for biharmonic
functions in Lipschitz and $C^1$ domains
\jour Commen. Math. Helv.
\vol 68 \yr 1993 \pages 385-414
\endref

\ref\key PV3
\by J.~Pipher and G.~Verchota
\paper
Dilation invariant estimates and the boundary Garding
inequality for higher order elliptic operators
\jour Ann. of Math. \yr 1995 \vol 142 \pages 1-38
\endref

\ref\key PV4
\by J.~Pipher and G.~Verchota
\paper
Maximum principle for the polyharmonic equation
on Lipschitz domains
\jour Potential Analysis\vol 4
\yr 1995
\pages 615-636
\endref

\ref \key S1
\by Z.~Shen
\paper Resolvent estimates in $L^p$ for elliptic systems
in Lipschitz domains
\jour J. Funct. Anal. \yr 1995 \vol 133 \pages 224-251
\endref

\ref \key S2
\by Z.~Shen
\paper A note on the Dirichlet problem
for the Stokes system in Lipschitz domains
\jour Proc. Amer. Math. Soc.
\vol 123 \yr 1995 \pages 801-811
\endref

\ref\key S3
\by Z.~Shen
\paper Boundary value problems in Morrey spaces
for elliptic systems on Lipschitz domains
\jour Amer. J. Math.
\pages 1079-1115
\vol 125\yr 2003
\endref

\ref\key S4
\by Z.~Shen
\paper
Weighted estimates for elliptic systems in Lipschitz domains
\jour to appear in Indiana Univ. Math. J.
\endref

\ref\key S5
\by Z.~Shen
\paper
Bounds of Riesz transforms on $L^p$ spaces for second order elliptic
operators
\jour to appear in Ann. Inst. Fourier (Grenoble)
\endref

\ref\key V1
\by G.~Verchota
\paper Layer potentials and regularity for the Dirichlet
problem for Laplace's equation
\jour J. Funct. Anal.
\vol 59 \yr 1984 \pages 572-611
\endref

\ref\key V2
\by G.~Verchota
\paper The Dirichlet problem for the biharmonic equation
in $C^1$ domains
\jour Indiana Univ. Math. J.
\vol 36 \yr 1987 \pages 867-895
\endref

\ref\key V3
\by G.~Verchota
\paper The Dirichlet problem for the polyharmonic equation
in Lipschitz domains
\jour Indiana Univ. Math. J.
\vol 39 \yr 1990 \pages 671-702
\endref

\endRefs
\enddocument

\end